\newfont{\Bbb}{msbm10 scaled\magstephalf}

\documentclass[11pt]{article}
\usepackage{amssymb, amsmath, indentfirst}
\usepackage{latexsym}
\usepackage{graphics,color}
\graphicspath{{D:/picture}}

\title{sdfgsdfg*********\\poopop**\thanks{sdlfj}}
\author{*******\\efdfsdf***\thanks{asdf}}
\date{**********\\easdft}
\renewcommand\dfrac{\displaystyle\frac}

\newcounter{theorem}[section]
\newcounter{corollary}[section]
\newcounter{lemma}[section]
\newcounter{proposition}[section]
\newcounter{definition}[section]
\newcounter{assumption}[section]

\newcounter{remark}

\newtheorem{theorem}{\indent Theorem}[section]
\newtheorem{definition}{\indent Definition }[section]
\newtheorem{lemma}{\indent Lemma}[section]
\newtheorem{proposition}{\indent Proposition}[section]
\newtheorem{corollary}{\indent Corollary}[section]

\begin{document}
\title{  \bf{Multivariate Lagrange Interpolation and
 an Application of Cayley-Bacharach Theorem for It}}
\author{Xue-Zhang Liang$^a$ \ Jie-Lin Zhang \ Ming Zhang \ Li-Hong Cui \\
 \small {\sl Institute of Mathematics, Jilin University,
Changchun, 130012, P. R. China}\\
 \small {\sl
$^a$liangxz@jlu.edu.cn}
}
\date{}
\maketitle

\footnotetext{ The paper supported by the NSF (60542002) of China and The 985 Program of Jilin University }

{\noindent {\bf Abstract }
}

In this paper,we deeply research Lagrange interpolation of $n$-variables and give an application of
Cayley-Bacharach theorem for it. We pose the concept of sufficient intersection about
$s (1\leq s\leq n)$ algebraic hypersurfaces in $n$-dimensional complex Euclidean
space and discuss the Lagrange interpolation along the algebraic
manifold of sufficient intersection. By means of some
theorems ( such as Bezout theorem, Macaulay theorem and so on ) we
prove the dimension for the polynomial space $P_m^{(n)}$ along  the
algebraic manifold $S=s(f_{1},\cdots,f_{s})$(where
$f_{1}(X)=0,\cdots,f_{s}(X)=0$ denote $s$ algebraic hypersurfaces
) of sufficient intersection and give a convenient expression
 for dimension calculation by using the backward
difference operator. According to Mysovskikh theorem, we give a proof
of the existence and a characterizing condition of properly posed set of nodes of arbitrary
degree for interpolation along an algebraic
manifold of sufficient intersection. Further we point out that for $s$ algebraic hypersurfaces
$f_{1}(X)=0,\cdots,f_{s}(X)=0$ of sufficient intersection,
the set of polynomials $f_{1},\cdots,f_{s}$ must constitute the
$H$-base of ideal $ I_{s}=<f_{1},\cdots,f_{s}>$. As a main result of this paper,
we deduce a general method of constructing properly
posed set of nodes for Lagrange interpolation along an algebraic
manifold, namely the superposition interpolation process. At the end of
the paper, we use the extended Cayley-Bacharach theorem to
resolve some problems of Lagrange interpolation along the
$0$-dimensional and $1$-dimensional algebraic manifold. Just the application
of Cayley-Bacharach theorem constitutes the start point of constructing properly
posed set of nodes along the high dimensional algebraic manifold by using
the superposition interpolation process.


\setcounter{equation}{0}
\setcounter{section}{1}
\thispagestyle{plain}
{\noindent {\bf 1. Introduction }
}

Multivariate Lagrange interpolation is an important problem of
computational mathematics and approximation theory. The research
of the theory and method of multivariate interpolation is
developing rapidly in the past few decades to solve some practical
problems such as the calculation of multivariate function, the
design of surface, the construction of finite element scheme, and
so on. In the research of multivariate Lagrange interpolation, the
well-posedness problem is the first problem. Thus there are two
ways for the research of multivariate Lagrange interpolation: one
way is to construct the properly posed set of nodes (or PPSN, for
short)for a given space of interpolating polynomials ( see
[1]$\sim$[8]). The other way is to find out the proper space of
interpolating polynomials for a given set of interpolation nodes,
specially to determine the interpolation
 space of minimal degree ( see [9]$\sim$[11]). The former is the main research aspect of this paper.

Interpolation is to solve the problem of constructing a function $p$ belonging to a finite dimensional linear space
to interpolate a given
set of data. Interpolation by univariate polynomial is a very classical topic. However, interpolation by  multivariate polynomials
is much more intricate and is an active area of research currently .

The thought of superposition interpolation can be ascended to the univariate Newton divided difference interpolation. To the bivariate case,
in 1948, Radon[12] gave the
Straight Line-Superposition Process for constructing the PPSN for bivariate polynomial space. Then Liang[5] deduced the
Superposition Interpolation Process for constructing the PPSN for $\mathbb{R}^{2}$( including Conic-Superposition Process )
in 1965. In 1977, Chung and Yao[1] founded the well-posedness of interpolation with the natural lattice in $\mathbb{R}^{n}$.
It is in 1981 that I.P.Mysovskikh[2] constructed Hyperplane-Superposition Process for constructing the PPSN for
$n$-variate polynomial space. M. Gasca[9] constructed the schemes of superposition interpolation for bivariate Lagrange
 and Hermite interpolation.
In 1991, L.Bos [13] put forward the superposition scheme of
constructing the PPSN for $n$-dimensional space $\mathbb{C}^{n}$
by using the PPSN along $(n-1)$-dimensional algebraic
hypersurface.  In 1998, Liang and l\"{u}[14] gave the method of
constructing the PPSN along a plane algebraic curve of arbitrary
degree by using the intersection of a straight line and the
algebraic curve . About in 2001,  Liang and l\"{u}[6] put forward
the Hyperplane-Superposition Process for constructing the PPSN for
Lagrange interpolation along the algebraic manifold in
$\mathbb{C}^{n}$. Carnicer and Gasca[8] also studied the problems
of bivariate Lagrange interpolation on conics and cubics in 2002.
In 2003, Bojanov and Xu discussed polynomial interpolation of two
variables based on points that are located on multiple circles in
[7]. Liang[20] deduced the Algebraic Curve-Superposition Process
for constructing the PPSN for interpolation along the plane
algebraic curve in $\mathbb {R}^{2}$ in 2004. Next year, Liang
researched the polynomial interpolation problems along the
algebraic surface and along the space algebraic curve in $\mathbb
{R}^{3}$ and gave Space Algebraic Curve-Superposition Process for
constructing the PPSN along the algebraic surface in $\mathbb
{R}^{3}$.

By generalizing the above results,  we extend the superposition interpolation process to the $n$-dimensional space $\mathbb{C}^{n}$ and
give some practical conclusions by using the extended Cayley-Bacharach theorem in this paper. Particularly, we deduce a recursive method of
constructing the PPSN for Lagrange interpolation along the algebraic manifold of arbitrary dimension in $\mathbb{C}^{n}$ (
including $\mathbb{C}^{n}$ ).

The content of this paper is following: The upper bound of
dimension for polynomial space along an algebraic manifold of
sufficient intersection is proved in \S2 and then we give an
convenient dimension expression by using the backward difference
operator. The existence of the PPSN of arbitrary degree for
Lagrange interpolation along the algebraic manifold of sufficient
intersection and the dimension of interpolation space are
discussed in \S3, in this section we also deduce a general method
of constructing the PPSN for Lagrange interpolation along an
algebraic manifold, namely the superposition interpolation
process. Then in \S4 we give the characterizing conditions of the
PPSN for Lagrange interpolation and the relations between it and
$H$-base. In \S5, we use the extended Cayley-Bacharach theorem to
resolve some problems of Lagrange interpolation along the
$0$-dimensional and $1$-dimensional algebraic manifold and find
out some practical results. We give some examples about the
applications of the superposition interpolation process and
Cayley-Bacharach theorem in \S6.

In this paper, we discuss the problem of Lagrange interpolation in complex space $\mathbb{C}^{n}$. $\mathbb{P}_m^{(n)}$ denotes
 the complex space of all $n$-variate polynomials of total degree $\leq m$, i.e.
\begin{equation*}
\mathbb{P}_m^{(n)}=\{\sum_{0\leq |\alpha|\leq m} c_{\alpha}
x_1^{\alpha_1}x_2^{\alpha_2}\cdots x_n^{\alpha_n}\quad | ~c_{\alpha}\in \mathbb{C}\}.
\end{equation*}
where $|\alpha|=\alpha_1+ \cdots +\alpha_n,\alpha_1,\cdots,\alpha_n$ denotes nonnegative integer.

Let $m$ be an arbitrary integer, $n$ be a nonnegative integer, and
$$\left (\begin{array}{lc}m\\ n\end{array}\right )=
 \left \{
 \begin{array}{lc}
0,~~~~~~~~~~~~~~when\  \ m<n;\\
\dfrac{m!}{n!(m-n)!},~when\  \ m\geq n.
\end{array}\right.$$
In this paper $e_m^{(n)}=\left ( \begin{array}{c} m+n \\ n \end{array}\right)$.

\begin{definition}
\label{def3.1}
~~ Let $\mathcal{A}$=
$\{Q^{(i)}\}^{e_m^{(n)}}_{i=1}$ be a set of $e_m^{(n)}$ distinct points in $\mathbb{C}^{n}$.  Given any set
$\{q_i \}^{e_m^{(n)}}_{i=1}$ of complex numbers, we
seek a polynomial $f(X)\in P^{(n)}_m$ (where $f(X)=f(x_1, x_2,\cdots,x_n
))$ satisfying
\begin{equation}\label{1}
f(Q^{(i)})=q_i,~~ i=1, \cdots, e_m^{(n)}.
\end{equation}
If for any given set $\{q_i\}^{e_m^{(n)}}_{i=1}$ of
complex numbers there always exists a unique solution for the
equation system (1.1), we call the Lagrange interpolation problem a
properly posed interpolation problem and call the corresponding
set $\mathcal{A}$=$\{Q^{(i)}\}^{e_m^{(n)}}_{i=1}$ of nodes a properly posed set of nodes ( or PPSN for short )
 for ${\bf P}^{(n)}_m$.
\end{definition}

In practical problems, the Lagrange interpolation along an algebraic hypersurface or an algebraic manifold be often considered. So
only to discuss the interpolation in $\mathbb{C}^{n}$ is not enough. Here is the definition of the algebraic hypersurface
and the algebraic manifold.

\begin{definition}
\label{def3.3}
Let $P = \{ p_1, \ldots, p_s\} \subset \mathbb{C}[x_1,\cdots,x_n]$, $s$ a positive integer. The algebraic variety
$$
S=s(p_1,\cdots,p_s) = \{ (a_1,\cdots,a_n) \in \mathbb{C}^n | ~p_i(a_1,\cdots,a_n)=0,
i=1,\ldots,s \}
$$
is a geometrical object and we call it an algebraic manifold defined by $p_1,\cdots,p_s$. If $s=1$, $S=s(p_1)$ is called
an algebraic hypersurface, and if the total degree of $p_1$ is 1 we call it an algebraic hyperplane.
\end{definition}

\vskip0.5truecm
\setcounter{equation}{0}
\setcounter{definition}{0}
\setcounter{theorem}{0}
\setcounter{section}{2}
\noindent {\bf 2. The upper bound of dimension for $\mathbb{P}_m^{(n)}$ along the algebraic manifold of sufficient intersection }

Because a general algebraic manifold is more complicated, we only discuss the case of sufficient intersection in this
paper.

\begin{definition}
\label{def3.1} Let $k_1, \cdots, k_s$ be nonnegative numbers. We
call that $s( 1\leq s\leq n )$ algebraic hypersurfaces
$f_{1}(X)=0, \cdots,f_{s}(X)=0$ of degree $k_1, \cdots, k_s$,
respectively, in $\mathbb C^n$ sufficiently intersect at an
algebraic manifold $S=s( f_{1}, \cdots, f_{s} )$, if there exist
$n-s$ algebraic hypersurfaces $f_{s+1}(X)=0, \cdots,f_{n}(X)=0$ of
degree $k_{s+1}, \cdots, k_n$, respectively, in $\mathbb C^n$ such
that these $n$ algebraic hypersurfaces meet exactly at $k_1 \cdots
k_n$ mutually distinct finite points in $\mathbb C^n$.
\end{definition}

The following extended Bezout Theorem will be used:

\begin{theorem}
\label{Bezoutth}
Let $f_1, \cdots, f_n$ be polynomials of degree $k_1, \cdots, k_n$, respectively, in $\mathbb C[ x_1, \cdots, x_n ]$. Then hypersurfaces
$f_{1}(X)=0, \cdots, f_{n}(X)=0$ meet exactly either at an infinite aggregate or at $k_1 \cdots k_n$ points.  The calculation
about the number of the latter points
includes the infinite point and the multiplicity of each point.
\end{theorem}

It is easy to prove the following proposition by using Bezout theorem:

\begin{proposition}
\label{Bezoutpro}
Let $f_{1}(X)=0,\cdots,f_{n}(X)=0$ be $n$ algebraic hypersurfaces of degree $k_1, \cdots, k_n$, respectively, in $\mathbb C^n$
of
sufficient intersection. Let $g_{i}(X)$ denote the sum of all monomials of degree $k_{i}$ of $f_{i}$ and call $g_{i}(X)$
the highest degree polynomial of $f_{i}$, $i=1, 2,\cdots, n$.
 Then the common zero points of $g_{1}(X), \cdots, g_{n}(X)$ are
 only at $( 0, \cdots, 0 )$.
\end{proposition}

\begin{definition}\label{dy1.18}
The manifold $\{a\in \mathbb{C}^n$$~|~f_i(a)=0,~i=1,\cdots,r\}$ of ideal $I=<f_1,\cdots,f_r>$ is defined as a set of
points in $\mathbb{C}^n$, namely the set of solutions for equation set $ \{f_i(X)=0,~i=1,\cdots,r\}$, and write it
as $U(f_1,\cdots,f_r)$. We call the dimension of the manifold of ideal $I=<f_1,\cdots,f_r>$ as the dimension of the ideal $I$, namely
$dim(I)=dim(U(I))$. If $d$ is the dimension of ideal $I$ in $\mathbb{C}^n$, then $n-d$ is called the rank of $I$.
\end{definition}

\begin{proposition}\label{46}[15]
If the rank of a homogeneous ideal $B=<g_1,\cdots,g_r>$ is $r$, then the rank of the ideal
$<g_1,\cdots,g_i> is \  \ $i$,~\forall i=1,\cdots,r-1$ .
\end{proposition}

\begin{theorem}
\label{the3.9}[16] ( Macaulay ) Let $\cal B$$=<g_1,\cdots,g_s>$ be
a homogeneous ideal, $s\leq n$, and the degree of $g_i$ be $k_i$.
If the rank of ideal $\cal B$ is $s$ and $d_m(s)$ denotes the
dimension for the linear space of the all homogeneous polynomials
of degree $m$ in $\cal B$, then
\begin{equation*}
\label{3.15}
d_m(s)=M(n-1,m)-h_m(s)
\end{equation*}
where $M(n-1,m)=\dfrac{(n-1+m)!}{(n-1)!m!}$, $h_m(s)$ is the coefficient of $t^m$ in the expansion of the function
$\psi_s(t)=(1-t^{k_1})\cdots(1-t^{k_s})(1-t)^{-n}$.
\end{theorem}

For all the $n$-variate monomials of degree $\leq m$, we permute
them according to the certain order of total degree ( by ascending
power order ) as $\varphi_t(X), t=1, 2,\cdots, e_{m}^{(n)}$.We
choose a maximal linearly independent subset of the set of
elementary item $\{X^{\alpha}g_i~|~|\alpha|+k_i=m,~i=1,\cdots,s\}$
of degree $m$ in $\mathcal B$ and label all the elements in the
subset. Obviously the number of the labelled elements is $d_m(s)$
according to the above theorem. Constructing a matrix $G_m$ whose
row vectors are the coefficients of the expansion in power of
degree $m$ for the labelled elementary items . The element located
at $( i,j )$ of $G_m$ is the coefficient of monomial
$\varphi_t(X)$ in the expansion of the $i$-th elementary item,
where $t=M( n, m-1 )+j$. We can see $G_m$ is a matrix of order
$d_m(s)\times M(n-1,m)$.

\begin{definition}
\label{def3.5} We select out $d_m(s)$ linear independent columns
from $G_m$and call the monomials of degree $m$ corresponding to
the $d_m(s)$ columns the selected monomials and write the indexing
set as $T_{m}(s)$. We call the monomials of degree $m$
corresponding to the remainder columns as the unselected monomials
and write the indexing set as $T_{m}^{'}(s)$.
\end{definition}

\begin{lemma}
\label{lemma3.1}[2]
Let $M_s=k_1+\cdots+k_s-n$. Then any selected monomial of degree $m$ can be denoted as the linear combination
of the unselected monomials of degree $m$ and the elementary items of degree $m$ in $\cal B$. Especially, the
selected monomial of degree $m$
must be the linear combination of elementary item of degree $m$ in $\cal B$ if $s=n$ and $m\geq M_n$.
\end{lemma}

We get the following theorem about the upper bound of dimension
for the polynomial space $\mathbb P_m^{(n)}$ along the algebraic
manifold of sufficient intersection.

\begin{theorem}
\label{the3.10} Let $s$ algebraic hypersurfaces without multiple
factors ( or AHWMF and AHWMFs for plural form, for short )
$f_1(X)=0, \cdots, f_s(X)=0$ of degree $k_1, \cdots, k_s$,
respectively, in $\mathbb C^n$ sufficiently intersect at the
algebraic manifold $S=s(f_1,\cdots,f_s)$. Let $g_i(X)$ be the
highest degree polynomial of $f_i(X)$, $i=1, \cdots, s$ and
$\mathcal A=<f_1,\cdots,f_s>$, $\mathcal B=<g_1,\cdots,g_s>$. Then
for any given polynomial $f(X)$ in $\mathbb P_m^{(n)}$, there
exist monomials $\varphi _j(X), j \in  \bigcup\limits_{t=0}^{m}
T_{t}^{'}(s)$,  such that $f(X)$ can be represented as follows:
\begin{equation*}
\begin{array}{ll} f(X)=\sum\limits_{t=0}^{m}[\sum\limits_{j\in
T'_{t}(s)}a_j^{(t)}\varphi_j(X)+\sum\limits_{j=1}^s
c_j^{(t)}(X)f_j(X)]
\end{array}
\end{equation*}
where $a_{j}^{(t)}\in  \mathbb C$ and $deg c_j^{(t)}(X)f_j(X)=t$.
And the upper bound of dimension for $ \mathbb P_m^{(n)}$ along
the algebraic manifold $S$ is $\sum\limits_{j=0}^m h_j(s)$, where
$h_j(s)$ denotes the coefficient of $t^j$ in the expansion of the
function $(1-t^{k_1})\cdots(1-t^{k_s})(1-t)^{-n}$.
\end{theorem}

{\bf Proof of Theorem \ref{the3.10}:}

First for any given polynomial $f(X)\in P_m^{(n)}$, we have the following expression:
\begin{equation}
\label{3.16}
f(X)=g^{(m)}(X)+g^{(m-1)}(X)+\cdots+g^{(1)}(X)+g^{(0)}(X)
\end{equation}
where $g^{(i)}(X)$ is a homogeneous polynomial of degree $i$, it is the linear combination of monomials of degree $i$
of $f(X)$, $i=0, 1, \cdots, m$.

Let $M_s=k_1+\cdots+k_s-n$, $L_s=\min\{k_1,\cdots,k_s\}$,
 $N=k_1\cdots k_n$.

According to Lemma 2.1, for given $s$ ( $1\leq s\leq n$ ), any homogeneous polynomial of degree $m$ can be presented as the linear combination of
the unselected monomials of degree $m$ and the elementary items of degree $m$. According to Theorem 2.2, the number of
elements in $T_{m}^{'}(s)$ is $h_{m}(s)$. We can express $g^{(m)}(X)$ as follows:

\begin{equation}
\label{3.21}\begin{array}{ll} g^{(m)}(X)=\sum\limits_{j\in
T'_{m}(s)}a_j^{(m)}\varphi_j(X)+\sum\limits_{j=1}^s
c_j^{(m)}(X)g_j(X)\end{array}
\end{equation}

where the second term denotes the linear combination of the elementary items of degree $m$ and
$degc_j^{(m)}(X)g_j(X)=m$.

Because $g_j(X)$ is the highest degree polynomial of $f_j(X)$, then

\begin{equation}
\label{3.22}
g_j(X)=f_j(X)-f_j^{(k_j-1)}(X)
\end{equation}
where $f_j^{(k_j-1)}(X)$ is defined as a polynomial of degree $< k_j$ of $f_j(X)$.
  Substituting ( \ref{3.22} )into ( \ref{3.21} ), and substituting the expression of $g^{(m)}(X)$
 into ( \ref{3.16} ), we get:

\begin{equation}
\label{3.24} f(X)=\sum \limits _{j\in
T'_{m}(s)}a_j^{(m)}\varphi_j(X)+\sum_{j=1}^s c_j^{(m)}(X)f_j(X)
+\tilde{f}(X).
\end{equation}
where $\tilde{f}$ denotes a remainder term of a polynomial of degree $\leq m-1$. By this way we reduce $f(X)$
to a sum of a polynomial of degree $m$ and a polynomial of degree $m-1$. Do the same way to $\tilde{f}(X)$ as it to
$f(X)$ by replacing $m$ by $m-1$.
We keep doing it and stop until the degree of remainder term is less than $L_{s}$. At this time all the monomials of degree
less than $L_s$ are the unselected monomials and the number of them is $M(n,L_s-1)=h_{L_s-1}(s)+\cdots+h_0(s)$. Then we
get an expression of $f(X)$:

\begin{equation}
\label{3.25}\begin{array}{ll} f(X)&=\sum\limits_{t=L_s}^{m}
[\sum\limits_{j\in
T'_{t}(s)}a_j^{(t)}\varphi_j(X)+\sum\limits_{j=1}^s
c_j^{(t)}(X)f_j(X)]
+\sum\limits_{j=1}^{M(n,L_s-1)}a_j\varphi_j(X)\\
&=\sum\limits_{t=0}^{m}[\sum\limits_{j\in
T'_{t}(s)}a_j^{(t)}\varphi_j(X)+\sum\limits_{j=1}^s
c_j^{(t)}(X)f_j(X)]
\end{array}
\end{equation}
where $\deg c_j^{(t)}(X)f_j(X)=t$ and the monomials of degree less than $L_s$ are unselected
monomial and the number of them is $M(n,L_s-1)=h_{l_{s}-1}(s)+\cdots+h_{0}(s)$.

Because the values of $f_i(X),~i=1,\cdots,s$ are zero along the
manifold $S=s(f_1,\cdots,f_s)$,  we see that the polynomial
$f(X)$of degree $m$ along the manifold can be expressed by the
linear combination of all the unselected monomials of degree $\leq
m$. The number of the unselected monomials is $h_0(s)+ \cdots
+h_m(s)$ and we write it as $H_m(s)=\sum\limits_{i=0}^m h_i(s)$.
So the upper bound of dimension for the polynomial space $\mathbb
P_m^{(n)}$ along  the algebraic manifold is
 $H_m(s)=h_0(s)+ \cdots+h_m(s)$.

Hence, we have the following conclusion: {\bf The upper bound of dimension for a polynomial space ${\mathbb P}_m^{(n)}$ along the
algebraic manifold $S=s(f_1,\cdots,f_s)$ is $H_m(s)=\sum\limits_{j=0}^{m} h_{j}(s)$}.

Thus we finish the proof of Theorem \ref{the3.10}.

We summarize some relations appeared in the above proof as follows:

\begin{eqnarray*}
\begin{array}{ll}
h_j(n)=0, & \hbox{when}\  \ ~j\geq M_n +1 \\
h_j(n)=M(n-1,j), & \hbox{when}\  \ ~0\leq j<L
\end{array}
\end{eqnarray*}
\begin{equation*}\label{3.30*}\begin{array}{ll}
H_m(n)=H_{M_{n}}(n)=\sum\limits_{j=0}^{M_{n}} h_j(n)=N,\\
\sum\limits_{j=0}^m d_j(n)=M(n,m)-N \end{array} \quad \hbox{when}\  \ ~m \geq M_n
\end{equation*}

To give an exact expression of the dimension for ${\mathbb P}_m^{(n)}$ along an algebraic  manifold of sufficient
intersection, we need the following $Mysovskikh$ Theorem:

\begin{theorem}\label{the3.13}[2](Mysovskikh)
Suppose $n$ algebraic hypersurfaces $f_{1}(X)=0,\cdots,f_{n}(X)=0$ of degree $k_1,\cdots,k_n$, respectively,
exactly meet at $N=k_1 \cdots k_n$ mutually distinct finite points $\{X^{(1)},\cdots,X^{(N)}\}$.
Let $I=\{ 1, 2, \cdots, N\}, J=\bigcup\limits_{m=0}^{M_{n}} T^{'}_{m}(n)$. Then there exist a set of $n$-variate monomials
of degree $m$:
\begin{equation*}
\{ \varphi_{j}(X) | j\in T^{'}_{m}(n), m=0, 1, \cdots, M_{n} \}
\end{equation*}
such that the $Vandermonde$ matrix:
\begin{equation}\label{max2.1}
(\varphi_{j}(X^{(i)}))_{i\in I,j\in J}
\end{equation}
is nonsingular.
\end{theorem}

We introduce the following backward difference operator and some
notations for convenience:

\begin{definition}
\label{def3.6}
Let $n,k_1,k_2,\cdots,k_{n},$ be natural numbers, we give the following notations:
\begin{equation*}
\begin{array}{lc}
e_m^{(n)}(k_1)=\nabla_{k_1} e_m^{(n)}=e_m^{(n)}-e_{m-k_1}^{(n)}\\
e_m^{(n)}(k_1,k_2)=\nabla_{k_2}\nabla_{k_1} e_m^{(n)}=\nabla_{k_2}(e_m^{(n)}-e_{m-k_1}^{(n)})=
\nabla_{k_2}e_m^{(n)}-\nabla_{k_2}e_{m-k_1}^{(n)}\\
\vdots \\
e_m^{(n)}(k_1,k_2,\cdots,k_{n})=\nabla_{k_{n}}\cdots\nabla_{k_2}\nabla_{k_1}e_m^{(n)}
\end{array}
\end{equation*}
\end{definition}

Then we have the following relations by calculation:
\begin{proposition}
\label{pro3.3} By calculation, we have:
\begin{equation*}
h_m(s)=e_m^{(n)}(k_1,k_2,\cdots,k_{s})-e_{m-1}^{(n)}(k_1,k_2,\cdots,k_{s})
\end{equation*}
And the upper bound of dimension for ${\mathbb P}_m^{(n)}$ along an algebraic manifold equals to $e_m^{(n)}(k_1,\cdots,k_s)$
given by Definition 2.4, i.e.:
\begin{equation*}\label{3.30}
H_m (s)=\sum\limits_{j=0}^m h_j(s)=e_m^{(n)}(k_1,\cdots,k_s)
\end{equation*}
\end{proposition}

\vskip0.5truecm
\setcounter{equation}{0}
\setcounter{definition}{0}
\setcounter{theorem}{0}
\setcounter {proposition}{0}
\setcounter{section}{3}
\noindent {\bf 3. The superposition interpolation, the existence of the PPSN for Lagrange interpolation and
             the dimension for the interpolation space ${\mathbb P}_m^{(n)}$ along the algebraic manifold }

\begin{definition}
\label{def3.6*}
Suppose $s (1\leq s\leq n)$ AHWMFs $f_1(X)=0, \cdots,f_s(X)=0$ of degree $k_1,\cdots,k_s$, respectively, sufficiently
intersect at the algebraic manifold $S=s(f_1,\cdots,f_s),$
Let $\mathcal A=\{Q^{(i)}\}_{i=1}^{e_m^{(n)}(k_1,\cdots,k_s)}$ be $e_m^{(n)}(k_1,\cdots,k_s)$ distinct points
on the algebraic manifold $S=s(f_1,\cdots,f_s)$. Given any set $\{q_i\}_{i=1}^{e_m^{(n)}(k_1,\cdots,k_s)}$ of complex
numbers, we are to seek a polynomial $f(X)\in \mathbb P_m^{(n)}$ satisfying:
\begin{equation}
\label{3.38}
f(Q^{(i)})=q_i\quad i=1,\cdots,e_m^{(n)}(k_1,\cdots,k_s)
\end{equation}

We call the set $\mathcal A=\{Q^{(i)}\}_{i=1}^{e_m^{(n)}(k_1,\cdots,k_s)}$ of nodes a PPSN for Lagrange interpolation
of degree $m$ along the algebraic manifold $S=s(f_1,\cdots,f_s)$ and write it as $\mathcal A\in I_m^{(n)}(S)$ ( where $I_m^{(n)}(S)$
denotes the set of all the PPSN for Lagrange interpolation of degree $m$ along the algebraic manifold $S$ ), if for each
given set $\{q_i\}_{i=1}^{e_m^{(n)}(k_1,\cdots,k_s)}$ of complex numbers there always exists a solution for the equation
system ( \ref{3.38} ).
\end{definition}

We can construct a PPSN for the $n$-dimensional space ${\mathbb P}_m^{(n)}$ by using the PPSN along the $(n-1)$-dimensional
hypersurface ( see [6]).

\begin{proposition}
Let $\mathcal A=\{Q^{(i)}\}_{i=1}^{e_m^{(n)}}$ be a PPSN for ${\mathbb P}_m^{(n)}$. Suppose a AHWMF $q(X)=0$ of degree $k$
does not pass through any points in $\mathcal A$. Choose arbitrarily a PPSN $\mathcal B$ of degree $m+k$ along the
hypersurface $q(X)=0$ such that $\mathcal B \in I^{(n)}_{m+k} (q)$. Then $\mathcal A \cup \mathcal B$ must be a PPSN
 for ${\mathbb P}_{m+k}^{(n)}$.
\end{proposition}

The following theorem is the superposition interpolation process to construct the PPSN for Lagrange interpolation along an algebraic
 manifold:

\begin{theorem}\label{the3.18}
Let $f_1(X)=0,\cdots,f_s(X)=0$ be $s$($1\leq s \leq n$) AHWMFs of
degree $k_1,\cdots,k_s$, respectively, in $\mathbb C^n$. Let
$e^{(n)}_{m}(k_{1},\cdots,k_{s})$ distinct points on the algebraic
manifold $S_{n-s}=s(f_1,\cdots,f_s)$ be a PPSN for Lagrange
interpolation of degree $m$ along the algebraic manifold $S_{n-s}$
and write the set of points as $E_{S_{n-s}}^{(m)}$; Let
$e^{(n)}_{m-k_{s}}(k_{1},\cdots,k_{s-1})$ distinct points on the
algebraic manifold $S_{n-s+1}=s(f_1,\cdots,f_{s-1})$ be a PPSN for
Lagrange interpolation of degree $m-k_{s}$ along the algebraic
manifold $S_{n-s+1}$ and write the set of points as
$E_{S_{n-s+1}}^{(m-k_s)}$. Suppose the surface $f_{s}(X)=0$ does
not contain any point in $E_{S_{n-s+1}}^{(m-k_s)}$. Then
$E_{S_{n-s}}^{(m)}\cup E_{S_{n-s+1}}^{(m-k_{s})}$ must be a PPSN
of degree $m$ for Lagrange interpolation along the algebraic
manifold $S_{n-s+1}$.
\end{theorem}

{\bf Proof of Theorem \ref{the3.18}: }

The number of points in $E_{S_{n-s}}^{(m)} \cup E_{S_{n-s+1}}^{(m-k_s)}$ is

\begin{equation*}
\begin{array}{ll}
e_m^{(n)}(k_1,\cdots,k_s)+e_{m-k_s}^{(n)}(k_1,\cdots,k_{s-1})=\nabla_{k_s}e_m^{(n)}(k_1,\cdots,k_{s-1})+e_{m-k_s}^{(n)}(k_1,\cdots,k_{s-1})\\
=e_m^{(n)}(k_1,\cdots,k_{s-1})-e_{m-k_s}^{(n)}(k_1,\cdots,k_{s-1})+e_{m-k_s}^{(n)}(k_1,\cdots,k_{s-1})
=e_m^{(n)}(k_1,\cdots,k_{s-1})
\end{array}
\end{equation*}
which is exactly equal to the number of points contained in a PPSN of degree $m$ along the algebraic manifold
$S_{n-s+1}=s(f_1,\cdots,f_{s-1})$ and write the $e_m^{(n)}(k_1,\cdots,k_{s-1})$ points as
$\{Q^{(i)}\}^{e_m^{(n)}(k_1,\cdots,k_{s-1})}_{i=1}$.

We will prove the theorem according to Definition \ref{def3.6*}.
For any given set of complex numbers
$\{q_i\}_{i=1}^{e_m^{(n)}(k_1,\cdots,k_{s-1})}$, we will seek a
polynomial $f(X)\in \mathbb P_m^{(n)}$ satisfying:
\begin{equation}\label{3.88}
f(Q^{(i)})=q_i\quad i=1,\cdots,e_m^{(n)}(k_1,\cdots,k_{s-1})
\end{equation}
Let $f(X)=g_m(X)+\alpha (X)f_s(X)\in \mathbb P_m^{(n)}$, where $g_m(X)$ denotes the interpolation polynomial of degree $m$
along the algebraic manifold $S_{n-s}$ and $\alpha (X)\in \mathbb P_{m-k_s}^{(n)}$.

We can get $g_m(X)$ and $\alpha (X)$ such that $f(X)$ satisfies the interpolation condition ( \ref{3.88} ).

Because $E_{S_{n-s}}^{(m)}$ is a PPSN of degree $m$ for Lagrange
interpolation along the algebraic manifold $S_{n-s}$, and
$f_s(Q^{(i)})=0, i=1,\cdots,e_m^{(n)}(k_1,\cdots,k_s)$, then:
\begin{equation*}
q_i=f(Q^{(i)})=g_m(Q^{(i)})\quad
i=1,\cdots,e_m^{(n)}(k_1,\cdots,k_s)
\end{equation*}
Due to Definition \ref{def3.6*}, we can get $g_m(X)$.

Because $E_{S_{n-s+1}}^{(m-k_s)}$ is a PPSN of degree $m-k_{s}$ for Lagrange interpolation along the algebraic manifold
$S_{n-s+1}$, then:
\begin{equation*}
q_i=f(Q^{(i)})=g_{m}(Q^{(i)})+\alpha (Q^{(i)})f_{s}(Q^{(i)})\quad
i=e_m^{(n)}(k_1,\cdots,k_s)+1,\cdots,e_{m}^{(n)}(k_1,\cdots,k_{s-1})
\end{equation*}

$f_s(Q^{(i)})\neq 0,g_m(X)$ is known by above process, so we get $\alpha (X)\in \mathbb P_{m-k_s}^{(n)}$.

And then, we get $f(X)\in \mathbb P_m^{(n)}$ satisfying the interpolation condition ( \ref{3.88} ). So $E_{S_{n-s}}^{(m)}\cup E_{S_{n-s+1}}^{(m-k_{s})}$
is a PPSN of degree $m$ along the algebraic manifold $S_{n-s+1}$ by Definition 3.1.

Thus we complete the proof of Theorem  \ref{the3.18}.

Next we give a series of proofs of the existence of a PPSN of arbitrary degree for Lagrange interpolation along an algebraic
manifold of sufficient intersection.

First we prove the existence of the PPSN of arbitrary degree for Lagrange interpolation along the $0$-dimensional
algebraic manifold
$S_0=s(f_1,\cdots,f_n)$.

\begin{theorem}\label{the3.14}
Suppose $n$ AHWMFs $f_1(X)=0, \cdots,f_n(X)=0$ of degree $k_1,\cdots,k_n$, respectively, in $\mathbb C^n$ sufficiently
intersect at $N=k_1\cdots k_n$ distinct points $\{Q^{(i)}\}_{i=1}^N$. Let $M_n=k_1+\cdots+k_n-n$. Then

(1) The set of these $N$ points $\{Q^{(i)}\}_{i=1}^N$ is a PPSN of degree $m$ for Lagrange interpolation along the
algebraic manifold $S_0=s(f_1,\cdots,f_n)$ when $m\geq M_n$.

(2) There must exist a PPSN of degree $m$ for Lagrange interpolation along the algebraic manifold $S_0=s(f_1,\cdots,f_n)$
contained in the set of these $N$ points $\{Q^{(i)}\}_{i=1}^N$ when $m<M_n$.
\end{theorem}

{\bf Proof of Theorem \ref{the3.14}: }

We will prove the theorem according to Definition  \ref{def3.6*}. For any given set of complex numbers
$\{q_i\}_{i=1}^{\nu}, \nu=\sum\limits_{j=0}^m h_j(n)$, we will seek a polynomial $f(X)\in \mathbb P_m^{(n)}$ satisfying:
\begin{equation}\label{3.45}
f(Q^{(i)})=q_i\quad i=1,\cdots,\nu
\end{equation}

(1) When $m\geq M_n$

Due to Theorem 2.2, we know $\nu=\sum\limits_{j=0}^{M_n}
h_j(n)=N$. Since for any given $f(X)\in \mathbb P_m^{(n)}$, we
have
\begin{equation}\label{3.46}
f(X)=\sum\limits_{j=1}^{e_{m}^{(n)}} c_{j}\varphi_{j}(X)
\end{equation}
the interpolation condition ( \ref {3.45} ) is an equation set of order $N$ and the coefficient matrix of it is a matrix of
order $N \times e_{m}^{(n)}$. The coefficient matrix must contain the matrix ( \ref{max2.1} ) as its submatrix of order $N \times N$.
According to Theorem 2.4, we know
the submatrix ( \ref{max2.1} ) is nonsingular. So the linear equation set (3.3)-(3.4) must have a solution, i.e. the interpolation
condition is satisfied.

(2) When $m<M_n$

At this time the interpolation condition ( \ref {3.45} ) is an
equation set of order $e_{m}^{(n)}(k_1,\cdots,k_n)$. Its
coefficient matrix is a matrix of order
$e_{m}^{(n)}(k_1,\cdots,k_n)\times e_{m}^{(n)}$. According to
Theorem 2.4, there must exist a sequence of nonsingular nested
submatrices $B_{m} ( m=M_n-1, \cdots,0 )$ of order
$e_{m}^{(n)}(k_1,\cdots,k_n) \times e_{m}^{(n)}(k_1,\cdots,k_n)$
contained in the previous $e_{m}^{(n)}(k_1,\cdots,k_n)$ columns of
the matrix ( \ref{max2.1} ). In this submatrix sequence, the
latter is the submatrix of the former. The set of
$e_{m}^{(n)}(k_1,\cdots,k_n)$ nodes corresponding to the rows of
submatrix $B_{m}$ constitutes a PPSN for interpolation of degree
$m$ along the algebraic manifold and write it as $E_{S_0}^{(m)}$.
The following relation is obvious:
\begin{equation*}
E_{S_0}^{(0)}\subset E_{S_0}^{(1)}\subset \cdots \subset
E_{S_0}^{(M_n-1)}\subset E_{S_0}^{(M_n)}=E_{S_0}^{(M_n+1)}\cdots
\end{equation*}
where $E_{S_0}^{(i)}$ denotes the PPSN for interpolation of degree $i$ along the algebraic manifold $S_0$.

Synthesizing the above results (1) and (2) we have that there exists a PPSN of arbitrary degree along the
$0$-dimensional algebraic manifold $S_0=s(f_1,\cdots,f_n)$ and the number of points contained in it is
$\sum\limits_{j=0}^{m} h_j(n)=e_{m}^{(n)}(k_1,\cdots,k_n)$

Thus we complete the proof of Theorem \ref{the3.14}.

Further we prove the existence of the PPSN of arbitrary degree for interpolation along the $1$-dimensional algebraic manifold
$S_1=s(f_1,\cdots,f_{n-1})$.

\begin{theorem}\label{the3.15}
Let $n-1$ AHWMFs $f_1(X)=0, \cdots, f_{n-1}(X)=0$ of degree $k_1, \cdots, k_{n-1}$, respectively, in
$\mathbb C^n$ sufficiently intersect at the algebraic manifold $S_1=s(f_1,\cdots,f_{n-1})$.
Let $M_{n-1}=k_1+\cdots+k_{n-1}-n$. Then there must exist a PPSN of arbitrary degree for Lagrange interpolation along this
algebraic manifold and the number of points contained in it is $\sum\limits_{j=0}^m h_j(n-1)$$=e_m^{(n)}(k_1,\cdots,k_{n-1})$.
Especially  $\sum\limits_{j=0}^m h_j(n-1)$ $=\dfrac{1}{2}k_1\cdots k_{n-1} (2m+n+1-k_1-\cdots -k_{n-1})$
when $m \geq M_{n-1}$.
\end{theorem}

{\bf Proof of Theorem \ref{the3.15}: }

Due to Definition \ref{def3.1}, we know there must exist an AHWMF $f_n(X)=0$ of degree $k_n$ which meets
$f_1(X)=0,\cdots,f_{n-1}(X)=0$ at $N=k_1\cdots k_n$ distinct points and write the set of nodes as $\cal A$.
From Theorem \ref{the3.14}, we know there exists a PPSN of arbitrary degree for interpolation along the $0$-dimensional algebraic manifold
$S_0=s(f_1,\cdots,f_n)$. Choose a PPSN of degree $k_n$ along $S_0$ and write it as $E_{S_0}^{(k_n)}=\{Q^{(i)}\}_{i=1}^{\nu}$.
Due to Theorem 3.2 and Proposition \ref{pro3.3}, we know the number of it is
$\nu=\sum\limits_{j=0}^{k_n}h_j(n)=e_{k_n}^{(n)}(k_1,\cdots,k_n)$. Choose arbitrarily a point $x^{(0)}$ in
$S_1=s(f_1,\cdots,f_{n-1})$ but not contained in $f_n(X)=0$. ( For example, we can choose $x^{(0)}$ as a solution of the
equation set $f_1(X)=0, f_2(X)=0\cdots,f_{n-1}(X)=0, f_n(X)+\epsilon=0$, where $\epsilon$ denotes a small positive real
number. ) Then it is obvious that the point $x^{(0)}$ is a PPSN
of degree $0$ for interpolation along $S_1$. From Theorem \ref{the3.18}, we see $\{x^{(0)}\}\cup E_{S_0}^{(k_n)}$ is a PPSN
of degree $k_n$ for interpolation along $S_1$ and write it as $E_{S_1}^{(k_n)}=E_{S_0}^{(k_n)}\cup\{x^{(0)}\}$. Then
$E_{S_1}^{(k_n)}\in I_{k_n}^{(n)}(S_1)$. The number of points contained in $E_{S_1}^{(k_n)}$ is $e^{(n)}_{k_{n}}(k_{1},\cdots,k_{n-1})$
and write it as $E_{S_1}^{(k_n)}=\{Q^{(i)}\}_{i=1}^{e^{(n)}_{k_{n}}(k_{1},\cdots,k_{n-1})}$. Because we can get
$e^{(n)}_{k_{n}}(k_{1},\cdots,k_{n-1})$ polynomials $\tilde{f}_{i}(X)$ for interpolation from the equation set :
\begin{eqnarray}\label{the88}
\tilde{f}_{i}(Q^{(j)})=\delta_{ij},\ \
 i,j=1,\cdots,e^{(n)}_{k_{n}}(k_{1},\cdots,k_{n-1})
\end{eqnarray}
According to (2.5), $\tilde{f}_{i}(X)$ can be represented as follows:
\begin{eqnarray}\label{the888}
\tilde{f}_{i}(X)=\sum\limits_{j\in
J_{m}(s)} \tilde{a}_{i,j}^{(m)}\varphi_j(X)+\sum\limits_{j=1}^s
\tilde{c}_{i,j}^{(m)}(X)f_j(X) \  \  \  \ J_m(s)=\bigcup\limits_{t=0}^{m}T^{'}_{t} (s)
\end{eqnarray}
Combining (\ref{the88}) with (\ref{the888}), we can see the matrix
\begin{eqnarray*}
(\varphi_{j}(Q^{(i)}))_{i,j=1,\cdots,e^{(n)}_{k_{n}}(k_{1},\cdots,k_{n-1}) }
\end{eqnarray*}
is nonsingular, where the $j$-th column of the matrix is
$(\varphi_{j}(Q^{(1)}),\cdots,\varphi_{j}(Q^{(e^{(n)}_{k_{n}}(k_{1},\cdots,k_{n-1}))}))^{\top}$.
 Choose a nonsingular submatrix of order $e^{(n)}_{k_{n}-1}(k_{1},\cdots,k_{n-1})$ in the previous
$e^{(n)}_{k_{n}-1}(k_{1},\cdots,k_{n-1})$ columns of the above matrix and the set of nodes corresponding to the
submatrix constitutes a PPSN for interpolation of degree $k_{n}-1$ along the algebraic manifold $S_{1}$. By this way
we can prove the existence of the PPSN of degree $k_{n}-2,\cdots,1$ for interpolation along $S_{1}$.
We can construct the PPSN of degree $k_{n}+1,k_{n}+2,\cdots$ along $S_{1}$ by using the superposition interpolation
process in Theorem 3.1 and based on the existence of the PPSN of degree $k_{n}, k_{n-1}, \cdots, 1$ along $S_{1}$.

Thus we complete the proof of Theorem \ref{the3.15}.

Generally, by the same method, we get the following theorem:

\begin{theorem}\label{the3.16}
Let $s ( 1\leq s\leq n )$ AHWMFs $f_1(X)=0, \cdots, f_s(X)=0$ of degree $k_1, \cdots, k_s$, respectively, in
$\mathbb C^n$ sufficiently intersect at algebraic manifold $S_{n-s}=s(f_1,\cdots,f_s)$. Let $m$ be a nonnegative integer.
Then there must exist a PPSN of degree $m$ for interpolation along $S_{n-s}$ ( which contains $e_m^{(n)}(k_1,
\cdots,k_s)=\sum\limits_{j=0}^m h_j(s)$ points and we write it as $E_{S_{n-s}}^{(m)}$ ).
\end{theorem}

Combining Theorem 3.4 with Theorem 2.3, we also get the following dimension theorem:

\begin{theorem}\label{the3.15*}
The dimension $dim~ ({\mathbb P}_m^{(n)}~|_{s(f_1,\cdots,f_s)}) ( where  1\leq s\leq n )$ for polynomial space
${\mathbb P}_m^{(n)}$ along the algebraic manifold
$S_{n-s}=s(f_1,\cdots,f_s)$ is expressed as follows:
$$
dim~({\mathbb P}_m^{(n)}~|_{s(f_1,\cdots,f_s)})=\sum\limits_{j=0}^{m} h_j(s)
$$
\end{theorem}

Thus, summarizing Theorem 3.1-3.5, we extend the superposition interpolation process to the algebraic manifold of
sufficient intersection and deduce a recursive method of constructing the PPSN of arbitrary degree for Lagrange
interpolation along the algebraic manifold of sufficient intersection in $\mathbb C^n$.

\vskip0.5truecm
\setcounter{equation}{0}
\setcounter{definition}{0}
\setcounter{theorem}{0}
\setcounter{section}{4}
\noindent {\bf 4. The property of the PPSN for Lagrange interpolation along algebraic manifold and $H$-base }

\begin{theorem}\label{the3.17}
Let $s( 1\leq s \leq n )$ AHWMFs $f_1(X)=0, \cdots, f_{s}(X)=0$ of degree $k_1, \cdots, k_{s}$, respectively, in
$\mathbb C^n$ sufficiently intersect at the algebraic manifold $S=s(f_1,\cdots,f_{s})$.
Suppose ${ E}_S^{(m)}=\{Q^{(i)}\}_{i=1}^{e_m^{(n)}(k_1,\cdots,k_s)}$ is a set of $e_m^{(n)}(k_1,\cdots,k_s)$
mutually distinct points
on the manifold $S=s(f_1,\cdots,f_s)$. Then ${ E}_S^{(m)}$ being a PPSN of degree $m$ for Lagrange interpolation along
$S=s(f_1,\cdots, f_s)$, if and only if, for any polynomial $g(X)\in \mathbb {\mathbb P}_m^{(n)}$ satisfying the
following zero-interpolation condition
\begin{equation*}
g(Q^{(i)})=0,\quad i=1,\cdots,e_m^{(n)}(k_1,\cdots,k_s)
\end{equation*}
there always exists the following decomposition:
\begin{equation*}
g(X)=\sum\limits_{i=1}^{s}\alpha_i(X)f_i(X).
\end{equation*}
where $\alpha_i(X)\in {\mathbb P}_{m-k_i}^{(n)}$ and $\alpha_i(X)\equiv 0$, $i=1,\cdots,s$ when $m<k_i$.
\end{theorem}

{\bf Proof of the Theorem \ref{the3.17}}

We will use the mathematical induction to prove the theorem.

(1) The result is obvious from [6] when $s=1$.

(2) Suppose the result is true when $s=j-1$. We will prove it being true when $s=j$.

$\Leftarrow$ Sufficiency. Let $E_{S_{n-j}}^{(m)}$$=\{Q^{(i)}\}_{i=1}^{e_m^{(n)}(k_1,\cdots,k_j)}$ be distinct points
on manifold $S_{n-j}=s(f_1,\cdots,f_j)$. We choose a PPSN of degree $m-k_j$ along
$S_{n-j+1}=s(f_1,\cdots,f_{j-1})$ which is beyond the hypersurface $f_j(X)=0$  and write it as
$E_{S_{n-j+1}}^{(m-k_j)}$$=\{Q^{(i)}\}_{i=e_m^{(n)}(k_1,\cdots,k_j)+1}^{e_m^{(n)}(k_1,\cdots,k_j)+e_{m-k_j}^{(n)}(k_1,\cdots,k_{j-1})}$.
Then we have $E_{S_{n-j+1}}^{(m-k_j)}$$\cap E_{S_{n-j}}^{(m)}=\emptyset$ and $E_{S_{n-j+1}}^{(m-k_j)}$$\in
I_{m-k_j}^{(n)}(S_{n-j+1})$. Due to Theorem 3.4, we know $E_{S_{n-j+1}}^{(m-k_j)}$ must exist. We can prove
$E_{S_{n-j}}^{(m)}\cup$$E_{S_{n-j+1}}^{(m-k_j)}$$\in I_m^{(n)}(S_{n-j+1})$ as follows:

From Definition \ref{def3.6*} and  \ref{def3.6} we know the number of points in
$E_{S_{n-j}}^{(m)}\cup$$E_{S_{n-j+1}}^{(m-k_j)}$ is
\begin{equation*}\begin{array}{ll}
e_m^{(n)}(k_1,\cdots,k_j)+e_{m-k_j}^{(n)}(k_1,\cdots,k_{j-1})=
\nabla_{k_j}e_m^{(n)}(k_1,\cdots,k_{j-1})+e_{m-k_j}^{(n)}(k_1,\cdots,k_{j-1})\\
=e_m^{(n)}(k_1,\cdots,k_{j-1})-e_{m-k_j}^{(n)}(k_1,\cdots,k_{j-1})+e_{m-k_j}^{(n)}(k_1,\cdots,k_{j-1})
=e_m^{(n)}(k_1,\cdots,k_{j-1})
\end{array}
\end{equation*}
which is exactly equal to the number of the points contained in a PPSN of degree $m$ along $S_{n-j+1}=s(f_1,\cdots,f_{j-1})$.

Suppose polynomial $g(X)\in {\mathbb P}_m^{(n)}$ satisfies zero-interpolation condition:
\begin{equation}\label{z3.61}
g(Q^{(i)})=0,~~~~\forall Q^{(i)}\in  E_{S_{n-j}}^{(m)}\cup
E_{S_{n-j+1}}^{(m-k_j)}
\end{equation}
Then from the hypothetical condition of the theorem we have:
\begin{equation}\label{z3.62}
g(X)=\sum\limits_{i=1}^j \alpha_i(X)f_i(X)
\end{equation}
where $\alpha_i(X)\in {\mathbb P}_{m-k_i}^{(n)}$ and $\alpha_i(X)\equiv 0$ when $m<k_i$, $i=1,\cdots,j$.

For $\forall Q^{(i)}$$\in E_{S_{n-j+1}}^{(m-k_j)}$, because $Q^{(i)}$ is not contained in $f_j(X)=0$, we get:
\begin{equation*}
g(Q^{(i)})=\alpha_j(Q^{(i)})f_j(Q^{(i)})=0,~~~\forall Q^{(i)}\in
E_{S_{n-j+1}}^{(m-k_j)}
\end{equation*}

Because $f_j(Q^{(i)})\neq 0$ when $\forall Q^{(i)}\in
E_{S_{n-j+1}}^{(m-k_j)}$, it is obvious that
$\alpha_j(Q^{(i)})=0$. From the founding of necessity for the
conclusion of inductive assumption for the case of $j-1$, we have
the following decomposition of $\alpha_j(X)$:
\begin{equation}\label{z3.63}
\alpha_j(X)=\sum\limits_{i=j}^{j-1} \beta_i(X)f_i(X)
\end{equation}
where $\beta_i(X)\in {\mathbb P}_{m-k_j-k_i}^{(n)}$ and
$\beta_i(X)\equiv 0$,$i=1,\cdots,j-1$ when $m<k_j+k_i$.

Substituting (\ref{z3.63}) into (\ref{z3.62}), we get
\begin{equation*}
\begin{array}{ll}
g(X)&=\alpha_j(X)f_j(X)+\sum\limits_{i=1}^{j-1} \alpha_i(X)f_i(X)=
\sum\limits_{i=1}^{j-1} \beta_i(X)f_i(X)\cdot f_j(X)+\sum\limits_{i=1}^{j-1} \alpha_i(X)f_i(X)\\
&=\sum\limits_{i=1}^{j-1}
(\beta_i(X)f_j(X)+\alpha_i(X))f_i(X)=\sum\limits_{i=1}^{j-1}
\gamma_i(X)f_i(X)
\end{array}
\end{equation*}
where $\gamma_i(X)=\beta_i(X)f_j(X)+\alpha_i(X)$ and $\gamma_i(X)\in {\mathbb P}_{m-k_i}^{(n)}$.
If $m<k_i$, then $\gamma_i(X)\equiv 0$, $i=1,\cdots,j-1$.

So from the inductive assumption for the case of $j-1$, we get
$E_{S_{n-j}}^{(m)}\cup E_{S_{n-j+1}}^{(m-k_j)}$ $\in
I_m^{(n)}(S_{n-j+1})$.

Thus due to Definition \ref{def3.6*}, we get: for an arbitrary set of complex numbers\\
$\{q_i\}_{i=1}^{e_m^{(n)}(k_1,\cdots,k_j)+e_{m-k_j}^{(n)}(k_1,\cdots,k_{j-1})}$, there always exists a polynomial
$f(X)\in {\mathbb P}_m^{(n)}$ satisfying:
\begin{equation*}
f(Q^{(i)})=q_i,~~~~i=1,\cdots,e_m^{(n)}(k_1,\cdots,k_j)+e_{m-k_j}^{(n)}(k_1,\cdots,k_{j-1})
\end{equation*}

In particular, the following is true to the subset $\{q_i\}_{i=1}^{e_m^{(n)}(k_1,\cdots,k_j)}$ :
\begin{equation*}
f(Q^{(i)})=q_i,~~~~i=1,\cdots,e_m^{(n)}(k_1,\cdots,k_j)
\end{equation*}

From Definition \ref{def3.6*}, we know $E_{S_{n-j}}^{(m)}$$=\{Q^{(i)}\}_{i=1}^{e_m^{(n)}(k_1,\cdots,k_j)}$
$\in I_m^{(n)}(S_{n-j})$. The sufficiency holds.

$\Rightarrow$ Necessity. The points of $E_{S_{n-j+1}}^{(m-k_j)}$ and $E_{S_{n-j}}^{(m)}$ are same as those in the
proof of the sufficiency, but the sets of $E_{S_{n-j+1}}^{(m-k_j)}$ and $E_{S_{n-j}}^{(m)}$ are both PPSN here.
 From Theorem 3.1, we can prove  $E_{S_{n-j}}^{(m)}\cup E_{S_{n-j+1}}^{(m-k_j)}$ $\in I_m^{(n)}(S_{n-j+1})$.

Because $E_{S_{n-j+1}}^{(m-k_j)}$$\in I_{m-k_j}^{(n)}(S_{n-j+1})$, then to the set of complex numbers \\
$\{g(Q^{(i)})/f_j(Q^{(i)})\}_{i=e_m^{(n)}(k_1,\cdots,k_j)+1}^{e_m^{(n)}(k_1,\cdots,k_j)+e_{m-k_j}^{(n)}(k_1,\cdots,k_{j-1})}$,
there always exists a polynomial $\alpha_j(X)\in {\mathbb P}_{m-k_j}^{(n)}$ satisfying:
\begin{equation*}
\alpha_j(Q^{(i)})=\dfrac{g(Q^{(i)})}{f_j(Q^{(i)})},~~~\forall
Q^{(i)} \in E_{S_{n-j+1}}^{(m-k_j)}
\end{equation*}
and $\alpha_j(X)\equiv 0$ when $m<k_j$.

Constructing a polynomial $q(X)\in {\mathbb P}_m^{(n)}$ as follows:
\begin{equation}\label{z3.68}
q(X)=g(X)-\alpha_j(X)f_j(X)
\end{equation}
It is obvious that $q(Q^{(i)})=0$ to $\forall Q^{(i)}$$\in E_{S_{n-j}}^{(m)}\cup E_{S_{n-j+1}}^{(m-k_j)}$.
Because $E_{S_{n-j}}^{(m)}\cup E_{S_{n-j+1}}^{(m-k_j)}$$\in I_m^{(n)}(S_{n-j+1})$ from the
conclusion of inductive assumption for the case of $j-1$, we know this $q(X)$
should have the following decomposition:

\begin{equation}\label{z3.69}
q(X)=\sum\limits_{i=1}^{j-1} \alpha_i(X)f_i(X)
\end{equation}
where $\alpha_i(X)\in {\mathbb P}_{m-k_i}^{(n)}$ and $\alpha_i(X)\equiv 0$ when $m<k_i$, $i=1,\cdots,j-1$.

Substituting (\ref{z3.69}) into (\ref{z3.68}), we get
\begin{equation*}
g(X)=q(X)+\alpha_j(X)f_j(X)=\sum\limits_{i=1}^{j-1}
\alpha_i(X)f_i(X)+\alpha_j(X)f_j(X)=\sum\limits_{i=1}^j
\alpha_i(X)f_i(X)
\end{equation*}
where $\alpha_i(X)\in {\mathbb P}_{m-k_i}^{(n)}$ and $\alpha_i(X)\equiv 0$ when $m<k_i$, $i=1,\cdots,j$.
This shows the validity of the case $s=j$.

Synthesizing the above proofs of (1) and (2) we know the conclusion is true to all $s(2\leq s \leq n)$.

Thus we complete the proof of Theorem \ref{the3.17}.

Next we give the relations between the algebraic hypersurfaces of sufficient intersection and the $H$-base of ideal.

\begin{definition}\label{dy1.12}[17]
We call a set of polynomials $G=\{g_1,\cdots,g_s\}$$\subset$$\mathbb{C}[x_1,\cdots,x_n]$$\setminus \{0\}$
a $H$-base of ideal $I=<g_1,\cdots,g_s>$ if for each polynomial $ p \in I$ of degree $m$, there always exist
$h_1,\cdots,h_s$$\in \mathbb{C}[x_1,\cdots,x_n]$ satisfying
\begin{equation}\label{1.1}
p=\sum\limits_{i=1}^s h_ig_i~~~\mbox{and}~deg(h_i)+deg(g_i)\leq deg(p),~~i=1,\cdots,s
\end{equation}
where the expression in (\ref{1.1}) is also called the $H$-expansion of $p$ about $G$.
\end{definition}

\begin{theorem}\label{the3.19}
Suppose $s( 1\leq s \leq n )$ AHWMFs $f_1(X)=0, \cdots, f_{s}(X)=0$ of degree $k_1, \cdots, k_{s}$, respectively, in
$\mathbb C^n$ sufficiently intersect at the algebraic manifold $S=s(f_1,\cdots,f_{s})$. Then the set of polynomials
$\{f_1,\cdots,f_s\}$ must be a $H$-base of ideal $I_s=<f_1,\cdots,f_s>$.
\end{theorem}

{\bf Proof of Theorem \ref{the3.19}: }

We will prove the theorem by using the inductive approach based on the Definition of $H$-base and Theorem 4.1.

(1) When $s=n$.

At this time the $n$ AHWMFs $f_1(X)=0, \cdots, f_{n}(X)=0$ of degree $k_1, \cdots, k_{n}$, respectively, in
$\mathbb C^n$ meet exactly at $N=k_1\cdots k_n$ distinct points. From Theorem \ref{the3.17}, we know that
to any polynomial $f(X)\in I_n$ of degree $m$, there must exist $\alpha_i(X)\in$${\mathbb P}_{m-k_i}^{(n)}$,
$i=1,\cdots,n$ satisfying:
\begin{equation}\label{03.51}
f(X)=\sum\limits_{i=1}^{n} \alpha_i(X)f_i(X)
\end{equation}

(2) When $s=n-1$.

At this time the $n-1$ AHWMFs $f_1(X)=0, \cdots, f_{n-1}(X)=0$ of degree $k_1, \cdots, k_{n-1}$, respectively, in
$\mathbb C^n$ just sufficiently intersect at the algebraic manifold $S_1=s(f_1,\cdots,f_{n-1})$. From the definition
of sufficient intersection, we know that there must exist an AHWMF $f_n(X)=0$ of degree $k_n$ which meet
$f_1(X)=0,\cdots,f_{n-1}(X)=0$ exactly at $N=k_1\cdots k_n$ distinct points and write the set of points as
 $\cal A$$=\{Q^{(i)}\}_{i=1}^{N}$.
Because $I_{n-1}=<f_1,\cdots,f_{n-1}>$ $\subset I_n=<f_1,\cdots,f_n>$, then for any given
polynomial $f(X)\in I_{n-1}\subset I_n$ of degree $m$ must have the expression as (\ref{03.51}) and
$deg   \alpha_i(X)f_i(X)\leq m$, $i=1,\cdots,n$.

Choose a PPSN
$E_{S_1}^{(m-k_n)}=$$\{Q^{(i)}\}_{i=1}^{e_{m-k_n}^{(n)}(k_1,\cdots,k_{n-1})}$
of degree $m-k_n$ along $S_{1}=s(f_1,\cdots,f_{n-1})$ beyond the
hypersurface $f_n(X)=0$ and writ it as $E_{S_1}^{(m-k_n)}\in
I^{(n)}_{m-k_{n}}(S_1)$. Then $E_{S_1}^{(m-k_n)}\cap$$\cal
A$$=\emptyset$. Substituting the points of $E_{S_1}^{(m-k_n)}$
into (4.7) we get

\begin{equation}\label{03.52}\begin{array}{llc}
f(Q^{(i)})&=\sum\limits_{j=1}^{n-1} \alpha_j(Q^{(i)})f_j(Q^{(i)})+\alpha_n(Q^{(i)})f_n(Q^{(i)})\\
&=\alpha_n(Q^{(i)})f_n(Q^{(i)})=0,~~~~\forall Q^{(i)}\in E_{S_1}^{(m-k_n)}.
\end{array}
\end{equation}
For $\forall Q^{(i)}\in E_{S_1}^{(m-k_n)}$, it is obvious that $f_n(Q^{(i)})\neq 0$, So $\alpha_n(Q^{(i)})=0$.
Due to $\alpha_n(X)\in$ ${\mathbb P}_{m-k_n}^{(n)}$ and Theorem \ref{the3.17} we know $\alpha_n(X)$ should have
the following decomposition:
\begin{equation}\label{03.53}
\alpha_n(X)=\sum\limits_{i=1}^{n-1} \beta_i(X)f_i(X)
\end{equation}
where $\beta_i(X)\in $${\mathbb P}_{m-k_n-k_i}^{(n)}$,
$i=1,\cdots,n-1$. Substituting (\ref{03.53}) into (\ref{03.52}),
we get
\begin{equation*}\begin{array}{llc}
f(X)&=\sum\limits_{i=1}^{n-1} \alpha_i(X)f_i(X)+\sum\limits_{i=1}^{n-1}\beta_i(X)f_i(X)\cdot f_n(X)\\
&=\sum\limits_{i=1}^{n-1}(\alpha_i(X)+\beta_i(X)f_n(X))f_i(X)=\sum\limits_{i=1}^{n-1}\gamma_i(X)f_i(X)
\end{array}
\end{equation*}
where $\gamma_i(X)=\alpha_i(X)+\beta_i(X)f_n(X)$ and $\gamma_i(X)\in $${\mathbb P}_{m-k_i}^{(n)}$,$i=1,\cdots,n-1$.
From Definition \ref{dy1.12}, we know $f_1,\cdots,f_{n-1}$ can constitute the $H$-base of ideal
$I_{n-1}=<f_1,\cdots,f_{n-1}>$.

(3)Using the same approach as (1) and (2) we can prove the case of $s=n-i,i=2,\cdots,n-1$ and we omit them here.

Synthesizing the above results of (1),(2) and (3) we complete the proof of the Theorem \ref{the3.19}.

\vskip0.5truecm \setcounter{equation}{0}
\setcounter{definition}{0} \setcounter{theorem}{0}
\setcounter{section}{5} \noindent {\bf 5.  The application of the
extended Cayley-Bacharach theorem in $\mathbb C^n$ to Lagrange
interpolation }

In this paper, we will use the extended Cayley-Bacharach theorem
in $\mathbb C^n$ to acquire some useful results. The
Cayley-Bacharach theorem is following:

\begin{theorem}
\label{1}[18]
Let $m,n$ and $r$ be nonnegative integer, $3\leq r\leq min\{m, n\}+2$. Suppose the algebraic curve $p(x_1, x_2)=0$
of degree $m$ and the algebraic curve $q(x_1, x_2)=0$ of degree $n$ meet exactly at $mn$ mutually distinct points and
write the set of points as $\cal A$. If a polynomial $f(x_1, x_2)\in \mathbb P^{(2)}_{m+n-2}$ such that $f(x_1, x_2)=0$
passes through the $mn-\dfrac{1}{2}(r-1)(r-2)$ points in $\cal A$, then it must pass through the $\dfrac{1}{2}(r-1)(r-2)$
remainder points, unless these $\dfrac{1}{2}(r-1)(r-2)$ remainder points lie on an algebraic curve of degree $r-3$.
\end{theorem}

The following is an extended Cayley-Bacharach theorem in $\mathbb C^n$ ( see [19] ).

\begin{theorem}\label{the3.20}
Let $k_1,\cdots,k_n$ be natural numbers and $n$ AHWMFs $f_1(X)=0,
\cdots, f_{n}(X)=0$ of degree $k_1, \cdots, k_{n}$, respectively,
in $\mathbb C^n$ meet exactly at $N=k_1\cdots k_n$ distinct points
and write the set of points as
 $\cal A$$=\{Q^{(i)}\}_{i=1}^N$.
Let $L=\min\{k_1,\cdots,k_n\}$, $M=k_1+\cdots+k_n-n$, $m$ be a natural number and $M-L+1\leq m\leq M-1$. If there exists
a polynomial $f(X)\in$${\mathbb P}_m^{(n)}$ of degree $m$  such that the surface $f(X)=0$ passes through the
$N-\dfrac{(n+M-m-1)!}{n!(M-m-1)!}$ points in $\cal A$, then it must pass through the $\dfrac{(n+M-m-1)!}{n!(M-m-1)!}$
remainder points, unless these $\dfrac{(n+M-m-1)!}{n!(M-m-1)!}$ remainder points lie on one algebraic hypersurface
of degree $M-m-1$.
\end{theorem}

The following theorem is the extended case of Theorem 5.2.

\begin{theorem}\label{the3.21*}
Let $n, m, k_1,\cdots,k_n$ be natural numbers,
$M=k_1+\cdots+k_n-n$, $0\leq m\leq M-1$. Suppose $n$ AHWMFs
$f_1(X)=0, \cdots, f_{n}(X)=0$ of degree $k_1, \cdots, k_{n}$,
respectively, in $\mathbb C^n$ meet exactly at $N=k_1\cdots k_n$
distinct points and write the set of points as
 $\cal A$$=\{Q^{(i)}\}_{i=1}^N$.  And
$\tilde{\cal A}$$\subset \cal A$ is a PPSN of degree $M-m-1$ along the $0$-dimensional algebraic manifold
$S_0=s(f_1,\cdots,f_n)$, $\tilde{\cal A}$ contains $e_{M-m-1}^{(n)}(k_1,\cdots,k_n)$ points all together,
$\cal B=$$\cal A \setminus \tilde{\cal A} $. If there exists a polynomial
$f(X)\in$${\mathbb P}_m^{(n)}$ of degree $m$ such that the hypersurface $f(X)=0$ passes through the
$N-e_{M-m-1}^{(n)}(k_1,\cdots,k_n)$ points in $\cal B$, then it must pass
through the $e_{M-m-1}^{(n)}(k_1,\cdots,k_n)$ remained points of $\cal A$.
\end{theorem}

By using the extended Cayley-Bacharach theorems for Lagrange interpolation along the $0$-dimensional
and $1$-dimensional algebraic manifolds, we get the following results:

\begin{theorem}\label{the3.21**}
Let $n, m, k_1,\cdots,k_n$ be natural numbers,
$M=k_1+\cdots+k_n-n$,$0\leq m\leq M-1$. Suppose $n$ AHWMFs
$f_1(X)=0, \cdots, f_{n}(X)=0$ of degree $k_1, \cdots, k_{n}$,
respectively, in $\mathbb C^n$ meet exactly at $N=k_1\cdots k_n$
distinct points and write the set of points as
 $\cal A$$=\{Q^{(i)}\}_{i=1}^N$. Suppose
$\tilde{\cal A}$$\subset \cal A$ is a PPSN of degree $M-m-1$ along the $0$-dimensional algebraic manifold
$S_0=s(f_1,\cdots,f_n)$. Then the remained points in $\cal B=$$\cal A \setminus \tilde{\cal A} $ must
be a PPSN for Lagrange interpolation of degree $m$ along the $0$-dimensional algebraic manifold $S_0=s(f_1,\cdots,f_n)$.
\end{theorem}

\begin{theorem}\label{the3.21}
Let $\{0\}={\bf P}^{(3)}_{-1}={\bf P}^{(3)}_{-2}=\cdots $, denote
the space of zero polynomial, and under these circumstances we
regard their corresponding PPSN as the empty set. Let
$k_1,\cdots,k_n$ be natural numbers and $L=\min
\{k_1,\cdots,k_n\}$. And let $m$ be an integer number satisfying
$m\leq L-1$. Suppose $n$ AHWMFs $f_1(X)=0, \cdots, f_{n}(X)=0$ of
degree $k_1, \cdots, k_{n}$, respectively, in $\mathbb C^n$ meet
exactly at $N=k_1\cdots k_n$ distinct points and write the set of
points as
 $\cal A$$=\{Q^{(i)}\}_{i=1}^N$.
$\mathcal{B}$ $\subseteq \mathcal{A}$ is a PPSN for ${\mathbb P}_{m}^{(n)}$. If ${\cal A}^t$ is a PPSN for Lagrange
interpolation of degree $M-m-k_t-1$ along the algebraic manifold $S_1^{(t)}=s(f_1,\cdots,f_{t-1},f_{t+1},\cdots,f_n)$
($1\leq t\leq n$) and ${\cal A}^t\cap \cal A$$=\emptyset$, then we have \\
$~~~~~~~~~~~~~~~~~~~~~~~~~~~~~~~~~~~~~~~~~~~~
{\cal A}^{t} \cup (\cal A$$ \setminus \cal B)$$\in I_{M-m-1}^{(n)}(S_1^{(t)})
$
\end{theorem}

From the above theorems, we deduce the following corollary which is convenient to use:
\begin{corollary}\label{3.10}
Let $k_1,\cdots,k_n$ be natural numbers, $m$ be an integer number
which satisfying $m<0$. Suppose $n$ AHWMFs $f_1(X)=0, \cdots,
f_{n}(X)=0$ of degree $k_1, \cdots, k_{n}$, respectively, in
$\mathbb C^n$ meet exactly at $N=k_1\cdots k_n$ distinct points
and write the set of points as
 $\cal A$$=\{Q^{(i)}\}_{i=1}^N$. If
 ${\cal A}^t$ is a PPSN for polynomial
interpolation of degree $M-m-k_t-1$ along the algebraic manifold $S_1^{(t)}=s(f_1,\cdots,f_{t-1},f_{t+1},\cdots,f_n)$
($1\leq t\leq n$) and ${\cal A}^t\cap \cal A$$=\emptyset$, then we have \\
$~~~~~~~~~~~~~~~~~~~~~~~~~~~~~~~~~~~~~~~~~~~~
{\cal A}^{t} \cup \cal A$$\in I_{M-m-1}^{(n)}(S_1^{(t)})
$
\end{corollary}

\vskip0.5truecm
\setcounter{equation}{0}
\setcounter{definition}{0}
\setcounter{theorem}{0}
\setcounter{section}{6}
\noindent {\bf 6.  Some examples}

In this section, we give some examples of the applications of the superposition interpolation process.

{\bf Example 1} Suppose $p(X)=0$ is a straight line in a give plane $\mathbb C^2$. Choose arbitrarily a point $Q^{(0)}$ in $p(X)=0$,
it is obvious $Q^{(0)}$ is a PPSN of degree $0$ along $p(X)=0$. Suppose a straight line $q_{1}(X)=0$ meets $p(X)=0$ exactly
at a point $Q^{(1)}$ such that $Q^{(0)}\neq Q^{(1)}$, then from Corollary 5.1 we know $\{Q^{(0)}\} \cup \{Q^{(1)}\}
\in I^{(2)}_{1}(p)$. By the same way, we can draw another straight line $q_{2}(X)=0$ meets $p(X)=0$ exactly at
a point $Q^{(2)}$ such that $Q^{(2)}\neq Q^{(0)}$ and $Q^{(2)}\neq Q^{(1)}$, from Corollary 5.1 we know
$\{Q^{(0)}, Q^{(1)}\} \cup \{Q^{(2)}\}\in I^{(2)}_{2}(p), \cdots$. Keeping doing this way, we can get: {\bf the $m+1$
mutually distinct points in the straight line $p(X)=0$ must be a PPSN of degree $m$ along $p(X)=0$ and write it as
$\{Q^{(0)}, Q^{(1)}, \cdots, Q^{(m)} \} \in I^{(2)}_{m}(p)$}( Example 1 (a)).  Further, using the superposition interpolation process
( Theorem 3.1 ) to this kind of PPSN along the above straight line,
we can construct a series of PPSN for $\mathbb P^{(2)}_m$ in $\mathbb C^2$, it is the Straight Line-Superposition Process
given by Radon ( Example 1 (b), see[12] ).

\begin{center}
\includegraphics*[0,0][500,300]{D: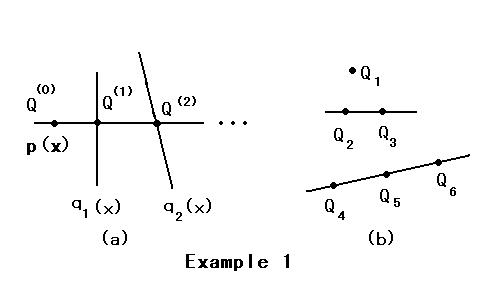}
\end{center}

{\bf Example 2} Suppose $p(X)=0$ is a conic in a given plane
$\mathbb C^2$. Choose arbitrarily a point $Q^{(0)}$ in $p(X)=0$,
it is obvious $Q^{(0)}$ is a PPSN of degree $0$ along $p(X)=0$.
Suppose a straight line $q_{1}(X)=0$ meets $p(X)=0$ exactly at two
points $Q^{(1)},Q^{(2)}$ such that $Q^{(1)}\neq Q^{(0)},
Q^{(2)}\neq Q^{(0)}$, then from Corollary 5.1 we know
$\{Q^{(0)}\}\cup \{Q^{(1)}, Q^{(2)}\} \in I^{(2)}_{1}(p)$. By the
same way, we can draw another straight line $q_{2}(X)=0$ meets
$p(X)=0$ exactly at two points $Q^{(3)}, Q^{(4)}$ such that
$Q^{(3)}, Q^{(4)}$ are distinct from the points $\{Q^{(0)},
Q^{(1)}, Q^{(2)}\}$, from Corollary 5.1 we know $\{Q^{(0)},
Q^{(1)}, Q^{(2)}\} \cup \{Q^{(3)}, Q^{(4)}\} \in I^{(2)}_{2}(p),
\cdots$. Keeping doing this way, we can get: {\bf the $2m+1$
mutually distinct points in the conic $p(X)=0$ must be a PPSN of
degree $m$ along $p(X)=0$ and write it as $\{Q^{(0)}, Q^{(1)},
\cdots, Q^{(2m)} \} \in I^{(2)}_{m}(p)$} ( Example 2 (a)).
Further, using the superposition interpolation process( Theorem
3.1 ) to this kind of PPSN along the above conic, we can construct
a series of PPSN for $\mathbb P^{(2)}_m$ in $\mathbb C^2$, it is
the Conic-Superposition Process in [5]. See Example 2 (b), here
$\{1, 2,\cdots, 6\}$ is a PPSN of degree $2$ for $\mathbb
P_2^{(2)}$, $\{1, 2,\cdots, 15\}$ is a PPSN of degree $4$ for
$\mathbb P_4^{(2)}$.

\begin{center}
\includegraphics*[10,15][460,275]{D: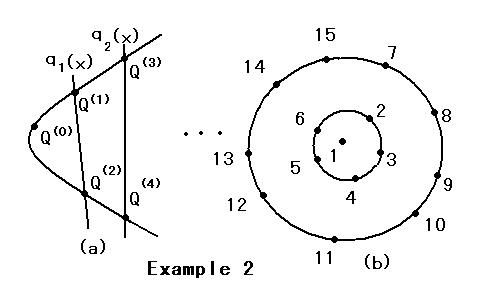}
\end{center}

{\bf Example 3} Suppose $3$ algebraic surface of degree $2\  \  \
p(X)=h_{1234}\bullet h_{5678}=0, q(X)=h_{1485}\bullet h_{2376}=0 $
and $r(X)=h_{1265}\bullet h_{4378}=0$ meet exactly at $8$ mutually
distinct points $\cal A$=$\{1, 2,\cdots, 8 \}$. It is obvious
$\cal A$ is not a PPSN for any subspace of $\mathbb P^{(3)}_{2}$.
From Theorem 5.4 we know $\{ 2,\cdots, 8\}$ must be a PPSN of
degree $2$ along the algebraic curve $C_{0}=s( p, q, r
)=l_{14}l_{23}l_{58}l_{67}
l_{15}l_{26}l_{48}l_{37}l_{12}l_{43}l_{56}l_{87}$. Choose
arbitrarily a point $1^{'}$ along the algebraic curve $C_{1}=s( p,
q )=l_{14}l_{23}l_{58}l_{67}$ such that $1^{'}$ is distinct from
any point in $\cal A$. From Theorem 5.5 we know $\{1^{'},
2,\cdots, 8\}\in I_{2}^{(3)}( C_1 )$. Then choose arbitrarily a
point $9$ on the plane $h_{1234}$ such that the point $9$ is
beyond the curve $C_1$. From Theorem 3.1, we know $\{1^{'},
2,\cdots, 8, 9\}\in I_{2}^{(3)}( p )$. Suppose the point $10$ is
any point in the space $\mathbb C^3$ beyond the plane $p(X)=0$.
Due to Theorem 3.1( or Proposition 3.1 ) we know the set of points
$\{1^{'}, 2,\cdots, 9, 10 \}$ must be a PPSN of degree $2$ for
$\mathbb P_2^{(3)}$.

\begin{center}
\includegraphics*[0,0][350,250]{D: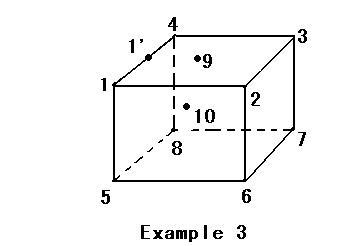}
\end{center}

Example 3 shows that by using the superposition interpolation process, we can construct a PPSN for the
algebraic surface from a series of PPSN along the algebraic curve, further we can construct a PPSN for the space $\mathbb C^n$.
Dimension increase is the substance of the superposition interpolation process. Similarly,
 the schemes for Lagrange interpolation such as given in [1], [2], [5],[6], [12]-[14],[20]
and so on can be deduced by the general superposition interpolation process put forward in this paper.
At the same time, it can be predicted that many new interpolation scheme in high dimensional space can also be deduced by the general
superposition interpolation process.

\newcounter{cankao}
\begin{list}
{[\arabic{cankao}]}{\usecounter{cankao}\itemsep=0cm}
\centerline{\bf References} \vspace*{0.5cm} \small \item K. C.
Chung, T. H. Yao. On lattices admitting unique Lagrange
interpolation. SIAM J. Numer, Anal., 1977,14:735$\sim$741.

\item I. P. Mysovskikh.
        Interpolatory cubature formulas (in
        Russian).Nauka, Moscow, 1981.

\item C. K. Chui, M.J.Lai.
        Vandermonde determinant and Lagrange interpolation in $\mathbb{R}^{s}$.
       in Nonlinear and Convex Analysis.B.L.Lin and S.Simon,eds.,Marcel Dekker,New York,
        1987, 23-36.

\item Y. Xu.
         Polynomial interpolation on the sphere and on the unit ball.
        Adv.in Comp. Math.
        2004 (20): 247-260.

\item X. Z. Liang.
       On the interpolation and approximation in several variables.
       Postgraduate Thesis, Jilin Univ.
       1965.

\item X. Z. Liang, C. M. L$\ddot{u}$, R. Z. Feng.
        Properly posed sets of nodes for multivariate Lagrange interpolation in $C^s$.
        SIAM, Numer. Anal.,
        2001,2(39): 578$\sim$595.

\item B.Bojanov. Y. Xu. On polynomial interpolation of two variable, J. Approx. Theory 2003, ( 120 ): 267-282.

\item J. M. Carnier, M. Gasca. Lagrange interpolation on conics
and cubics. CAGD 2002 ( 19 ): 313-326.

\item M. Gasca, J. I. Maeztu. On Lagrange and Hermite in $\mathbb{R}^{k}$. Numer. Math.1982(39):1-14.

\item  C. de. Boor, A. Ron.  On multivariate polynomial
interpolation.
        Constructive Approximation.   1990,6: 287$\sim$302.

\item T. Sauer. Polynomial interpolation of minimal degree and Gr$\ddot{o}$bner bases.
          Cambridge University Press,London Math.Soc. Lecture Notes,
          1998,251:483$\sim$494.

\item J. Radon. Zur mechanischen kubatur. Monatsh.~Math.,
       1948,52(4):286$\sim$300.

\item L. Bos.   On certain configurations of points in
$\mathbb{R}^n$ which are unisolvent for polynomial interpolation.
        SIAM J.Approx. Theory ,
        1991,64:271$\sim$280.

\item X. Z. Liang, C. M. L$\ddot{u}$.
       Properly posed set of nodes for bivariate interpolation.
      Approximation Theory IX,
      Computation Aspect,Vanderbit University Press, 1998,2:189$\sim$196.

\item B. Renschuch.
      Elementare and praktische idealtheorie.
      Berlin: VEB Deutscher Verlag der Wissenschaften.
     1976

\item F. S. Macaulay.
       Algebraic theory of modular systems. Cambridge tracts in mathematics and mathematical
        physics, Cambridge University press.
        1916, 19.

\item C. David, John Little, Donal O'sha.  Ideals,Varieties and Algorithms.
        Springer-Verlag, New York,
            1992.
\item J. G. Semple, L. Roth. Introduction to Algebraic Geometry. Oxford at the Clarendon Press, 1949.

\item D. Eisenbud, M. Green and J. Harris. Cayley-Bacharach theorems and conjectures.
       Bulletin of the American Mathematical Society.
      1996,Vol.33,Number 3:295-324.

\item X. Z. Liang, L. H. Cui, J. L. Zhang. The application of Cayley-Bacharach theorem to bivariate Lagrange interpolation.
       Journal of Computational and Applied Mathematics. 2004, 163(1):177-187.

\end{list}

\end{document}